\documentclass[11pt]{article}
\usepackage{latexsym}
\usepackage{amsfonts}
\usepackage{graphicx}
\usepackage{epstopdf}
\usepackage{color}
\usepackage[sort]{natbib}
\usepackage{lineno}

\newtheorem{theorem}{Theorem}
\newtheorem{lemma}{Lemma}
\newtheorem{corol}{Corollary}
\newtheorem{property}{Property}

\newcommand{\be}{\begin{equation}}
\newcommand{\ee}{\end{equation}}

\newcommand{\narrow}{\setlength{\itemsep}{1pt}
	\setlength{\parskip}{0pt}\setlength{\parsep}{0pt}}
\usepackage[normalem]{ulem}
\begin{document}
	\title{\bf The Obnoxious Facilities Planar $p$-Median Problem}
	\author{Pawel Kalczynski and Zvi Drezner\\
		Steven G. Mihaylo College of Business and Economics\\
		California State University-Fullerton\\
		Fullerton, CA 92834.\\e-mail: pkalczynski@fullerton.edu; zdrezner@fullerton.edu
	}
	\date{}
	\maketitle


	
\begin{abstract}
	In this paper we propose the planar obnoxious facilities $p$-median problem. In the $p$-median problem the objective is to find $p$ locations for facilities that minimize the weighted sum of distances between demand points and their closest facility. In the obnoxious version we add constraints that each facility must be located at least a certain distance from a partial set of demand points because they generate nuisance affecting these demand points. The resulting problem is extremely non-convex and traditional non-linear solvers such as SNOPT are not efficient. An efficient solution method based on Voronoi diagrams is proposed and tested. We also constructed the efficient frontiers of the test problems to assist the planers in making location decisions.
\end{abstract}
\noindent{\it Key Words: Location; Obnoxious Facilities; Continuous Location; Voronoi Diagrams.}

	\renewcommand{\baselinestretch}{1.6}
	\renewcommand{\arraystretch}{0.625}
	\large
	\normalsize
	
\section{Introduction}
	
Obnoxious location problems involve locating one or more facilities as far as possible from a set of communities \citep{GD75a,GD75b,Sham}. Most papers investigate the problem on networks or in discrete space \citep{CG78,DDS09,Erk:89,Welch1,CM79,CGMD16,Bat}; location in the interior of a network \citep{DW96,DDS09}; location on the plane \citep{WSD06,Sham,HPT81,DrWes83b}; location on the sphere \citep{Pap23}. Applications may include nuisance generated by the facilities such as airports, pollution generating industrial facilities, prisons, and other affecting residents living in a set of communities. Another type of applications assume that the nuisance is generated by the communities and the facilities should be located at locations with minimum nuisance. For example, the location of  schools or  hospitals which require a low noise level caused by a set of points or locating a telescope as far as possible from light sources. In most of these applications the nuisance propagates ``by air" and not along network links and thus are best modeled using  Euclidean distances.

While facilities may cause nuisance to demand points, commercial facilities aim to provide the best service possible. The commonly used objective is minimizing the weighted sum of distances between demand points to their closest facility. This model is termed the $p$-median model \citep{Das95,KH79med,DM15} also called the multi-source Weber problem  \citep{BHMT00,KS72}. Recent references for the planar $p$-median problem are \citep{DBMS13,DS15}. Recently, \citet{DDS18a} proposed the single facility median problem, the Weber problem \citep{W09,DKSW02,W93,LMW88}, with a minimum distance $D$ requirement between the facility and demand points. In this paper we investigate the multiple facility model in the plane that requires a minimum distance $D$ between some demand points and facilities.

In many applications demand points are getting services from the facilities but the facilities are also ``obnoxious" which means they also create nuisance to some or all demand points. This is addressed in many papers by having negative weights \citep[for example, ][]{DW91,MF93,TP89} in the Weber location problem \citep{W09}. In this paper we propose to treat obnoxious facilities by requiring a minimum distance between obnoxious facilities and some demand points. Two sets of points exist in the plane.  The set $N_d$ includes points that generate demand $w_i>0$ for $i\in N_d$, and the set $N_o$ includes points that the facility is obnoxious to them. The set of all points is $N=N_d\cup N_o$. It is possible that $N_d=N_o$ when all points generate demand and are also negatively affected by the facilities, or $N_d\cap N_o=\emptyset$, and all possibilities in between. The objective is to minimize the total distance traveled between demand points in $N_d$ and their closest facility while locating the facilities at least a distance $D$ from points in $N_o$. This is a better approach than locating facilities with attraction and repulsion. Let say that a demand point generates demand and is also negatively affected by the facilities. Should the weight of such a demand point be reduced by some negative value? How do we quantify the weight reduction? It is more realistic to minimize the total cost (generated by points in $N_d$) to satisfy the demand but require a minimum distance between facilities and demand points in $N_o$ that are negatively affected by the facilities.
 Many other objectives such as squared Euclidean distance, $p$-center, competitive facilities locations can yield useful models by adding the minimum distance requirement.

A related field of research is location among forbidden regions. In forbidden regions models, a facility is not allowed to be located in some regions in the plane. Demand points may be located there and travel is allowed through such regions. Therefore, the distance measure is not altered.
Papers that consider forbidden regions solve problems with convex objectives \citep{HS97,AnPar,BatGhose,ButCav,KatCoo}.
The optimal solution without forbidden regions can be found. If it is not in a forbidden region, this is the optimal solution. For convex problems it is shown that if the unconstrained solution is in a forbidden region, then the optimal solution is on the boundary of a forbidden region. This property does not necessarily hold for non-convex objectives. For example, the second best local optimum for the unconstrained problem may be feasible and better than all solutions on the boundaries of the forbidden regions.

The paper is organized as follows. In the next section we investigate the sensitivity of the solution to the value of the minimum distance $D$. In Section \ref{sec3} we show how to identify the feasible domain. In Section \ref{sec4} we detail the proposed solution algorithm. In Section \ref{comp} we test the solution algorithm and compare it to standard non-linear solvers on a set of test problems. In Section \ref{Efficient} we find and depict the efficient frontier. We conclude the paper and propose ideas for future research in Section \ref{concl}.

\section{The Effect of the value of $D$ on the Solution Approach}

	\begin{figure}[htp!]
	\begin{center}
		\setlength{\unitlength}{0.9in}
						\begin{picture}(1.1,1.15)
		\thicklines
		\put(0,0){\line(1,0){1}}
		\put(0,0){\line(0,1){1}}
		\put(0,1){\line(1,0){1}}
		\put(1,0){\line(0,1){1}}
		\put(0.00097,0.00367){\circle*{0.022}}
		\put(0.85243,0.84373){\circle*{0.022}}
		\put(0.84217,0.53687){\circle*{0.022}}
		\put(0.47523,0.01453){\circle*{0.022}}
		\put(0.83537,0.54207){\circle*{0.022}}
		\put(0.38603,0.55333){\circle*{0.022}}
		\put(0.90057,0.13927){\circle*{0.022}}
		\put(0.06483,0.74013){\circle*{0.022}}
		\put(0.15777,0.64847){\circle*{0.022}}
		\put(0.79163,0.65493){\circle*{0.022}}
		\put(0.92697,0.58967){\circle*{0.022}}
		\put(0.64643,0.17773){\circle*{0.022}}
		\put(0.72817,0.68287){\circle*{0.022}}
		\put(0.50923,0.98853){\circle*{0.022}}
		\put(0.28137,0.84807){\circle*{0.022}}
		\put(0.06003,0.56733){\circle*{0.022}}
		\put(0.50657,0.20527){\circle*{0.022}}
		\put(0.77883,0.19413){\circle*{0.022}}
		\put(0.52377,0.07447){\circle*{0.022}}
		\put(0.94563,0.94893){\circle*{0.022}}
		\put(0.65297,0.97567){\circle*{0.022}}
		\put(0.64043,0.71173){\circle*{0.022}}
		\put(0.41417,0.62887){\circle*{0.022}}
		\put(0.74323,0.16253){\circle*{0.022}}
		\put(0.52737,0.95407){\circle*{0.022}}
		\put(0.93403,0.78133){\circle*{0.022}}
		\put(0.91257,0.07127){\circle*{0.022}}
		\put(0.69283,0.84813){\circle*{0.022}}
		\put(0.68977,0.30047){\circle*{0.022}}
		\put(0.29963,0.44293){\circle*{0.022}}
		\put(0.17897,0.16167){\circle*{0.022}}
		\put(0.83443,0.44573){\circle*{0.022}}
		\put(0.90017,0.37487){\circle*{0.022}}
		\put(0.17723,0.53653){\circle*{0.022}}
		\put(0.57337,0.86007){\circle*{0.022}}
		\put(0.00803,0.19533){\circle*{0.022}}
		\put(0.11857,0.73727){\circle*{0.022}}
		\put(0.80683,0.70213){\circle*{0.022}}
		\put(0.65577,0.32647){\circle*{0.022}}
		\put(0.85363,0.13693){\circle*{0.022}}
		\put(0.50497,0.14767){\circle*{0.022}}
		\put(0.22843,0.37973){\circle*{0.022}}
		\put(0.18617,0.92087){\circle*{0.022}}
		\put(0.81123,0.11053){\circle*{0.022}}
		\put(0.41937,0.56607){\circle*{0.022}}
		\put(0.28203,0.80933){\circle*{0.022}}
		\put(0.12457,0.20327){\circle*{0.022}}
		\put(0.12083,0.75613){\circle*{0.022}}
		\put(0.42177,0.15247){\circle*{0.022}}
		\put(0.60763,0.03093){\circle*{0.022}}
		\put(0.63097,0.93367){\circle*{0.022}}
		\put(0.82243,0.51373){\circle*{0.022}}
		\put(0.27217,0.26687){\circle*{0.022}}
		\put(0.64523,0.88453){\circle*{0.022}}
		\put(0.06537,0.07207){\circle*{0.022}}
		\put(0.75603,0.62333){\circle*{0.022}}
		\put(0.93057,0.46927){\circle*{0.022}}
		\put(0.63483,0.01013){\circle*{0.022}}
		\put(0.98777,0.77847){\circle*{0.022}}
		\put(0.56163,0.12493){\circle*{0.022}}
		\put(0.55697,0.51967){\circle*{0.022}}
		\put(0.61643,0.84773){\circle*{0.022}}
		\put(0.15817,0.41287){\circle*{0.022}}
		\put(0.67923,0.85853){\circle*{0.022}}
		\put(0.51137,0.37807){\circle*{0.022}}
		\put(0.43003,0.63733){\circle*{0.022}}
		\put(0.53657,0.53527){\circle*{0.022}}
		\put(0.34883,0.46413){\circle*{0.022}}
		\put(0.35377,0.20447){\circle*{0.022}}
		\put(0.71563,0.41893){\circle*{0.022}}
		\put(0.28297,0.90567){\circle*{0.022}}
		\put(0.61043,0.38173){\circle*{0.022}}
		\put(0.84417,0.35887){\circle*{0.022}}
		\put(0.91323,0.03253){\circle*{0.022}}
		\put(0.75737,0.48407){\circle*{0.022}}
		\put(0.30403,0.85133){\circle*{0.022}}
		\put(0.94257,0.40127){\circle*{0.022}}
		\put(0.26283,0.11813){\circle*{0.022}}
		\put(0.51977,0.43047){\circle*{0.022}}
		\put(0.06963,0.91293){\circle*{0.022}}
		\put(0.80897,0.09167){\circle*{0.022}}
		\put(0.80443,0.11573){\circle*{0.022}}
		\put(0.33017,0.10487){\circle*{0.022}}
		\put(0.34723,0.40653){\circle*{0.022}}
		\put(0.80337,0.39007){\circle*{0.022}}
		\put(0.37803,0.26533){\circle*{0.022}}
		\put(0.14857,0.06727){\circle*{0.022}}
		\put(0.37683,0.97213){\circle*{0.022}}
		\put(0.48577,0.45647){\circle*{0.022}}
		\put(0.62363,0.60693){\circle*{0.022}}
		\put(0.13497,0.07767){\circle*{0.022}}
		\put(0.19843,0.04973){\circle*{0.022}}
		\put(0.61617,0.65087){\circle*{0.022}}
		\put(0.98123,0.98053){\circle*{0.022}}
		\put(0.64937,0.09607){\circle*{0.022}}
		\put(0.65203,0.87933){\circle*{0.022}}
		\put(0.15457,0.53327){\circle*{0.022}}
		\put(0.69083,0.02613){\circle*{0.022}}
		\put(0.25177,0.28247){\circle*{0.022}}
		\put(0.37763,0.50093){\circle*{0.022}}
		\end{picture}
		\begin{picture}(1.1,1.15)
		\thicklines
		\put(0,0){\line(1,0){1}}
		\put(0,0){\line(0,1){1}}
		\put(0,1){\line(1,0){1}}
		\put(1,0){\line(0,1){1}}
		\put(0.00097,0.00367){\circle*{0.075}}
		\put(0.85243,0.84373){\circle*{0.075}}
		\put(0.84217,0.53687){\circle*{0.075}}
		\put(0.47523,0.01453){\circle*{0.075}}
		\put(0.83537,0.54207){\circle*{0.075}}
		\put(0.38603,0.55333){\circle*{0.075}}
		\put(0.90057,0.13927){\circle*{0.075}}
		\put(0.06483,0.74013){\circle*{0.075}}
		\put(0.15777,0.64847){\circle*{0.075}}
		\put(0.79163,0.65493){\circle*{0.075}}
		\put(0.92697,0.58967){\circle*{0.075}}
		\put(0.64643,0.17773){\circle*{0.075}}
		\put(0.72817,0.68287){\circle*{0.075}}
		\put(0.50923,0.98853){\circle*{0.075}}
		\put(0.28137,0.84807){\circle*{0.075}}
		\put(0.06003,0.56733){\circle*{0.075}}
		\put(0.50657,0.20527){\circle*{0.075}}
		\put(0.77883,0.19413){\circle*{0.075}}
		\put(0.52377,0.07447){\circle*{0.075}}
		\put(0.94563,0.94893){\circle*{0.075}}
		\put(0.65297,0.97567){\circle*{0.075}}
		\put(0.64043,0.71173){\circle*{0.075}}
		\put(0.41417,0.62887){\circle*{0.075}}
		\put(0.74323,0.16253){\circle*{0.075}}
		\put(0.52737,0.95407){\circle*{0.075}}
		\put(0.93403,0.78133){\circle*{0.075}}
		\put(0.91257,0.07127){\circle*{0.075}}
		\put(0.69283,0.84813){\circle*{0.075}}
		\put(0.68977,0.30047){\circle*{0.075}}
		\put(0.29963,0.44293){\circle*{0.075}}
		\put(0.17897,0.16167){\circle*{0.075}}
		\put(0.83443,0.44573){\circle*{0.075}}
		\put(0.90017,0.37487){\circle*{0.075}}
		\put(0.17723,0.53653){\circle*{0.075}}
		\put(0.57337,0.86007){\circle*{0.075}}
		\put(0.00803,0.19533){\circle*{0.075}}
		\put(0.11857,0.73727){\circle*{0.075}}
		\put(0.80683,0.70213){\circle*{0.075}}
		\put(0.65577,0.32647){\circle*{0.075}}
		\put(0.85363,0.13693){\circle*{0.075}}
		\put(0.50497,0.14767){\circle*{0.075}}
		\put(0.22843,0.37973){\circle*{0.075}}
		\put(0.18617,0.92087){\circle*{0.075}}
		\put(0.81123,0.11053){\circle*{0.075}}
		\put(0.41937,0.56607){\circle*{0.075}}
		\put(0.28203,0.80933){\circle*{0.075}}
		\put(0.12457,0.20327){\circle*{0.075}}
		\put(0.12083,0.75613){\circle*{0.075}}
		\put(0.42177,0.15247){\circle*{0.075}}
		\put(0.60763,0.03093){\circle*{0.075}}
		\put(0.63097,0.93367){\circle*{0.075}}
		\put(0.82243,0.51373){\circle*{0.075}}
		\put(0.27217,0.26687){\circle*{0.075}}
		\put(0.64523,0.88453){\circle*{0.075}}
		\put(0.06537,0.07207){\circle*{0.075}}
		\put(0.75603,0.62333){\circle*{0.075}}
		\put(0.93057,0.46927){\circle*{0.075}}
		\put(0.63483,0.01013){\circle*{0.075}}
		\put(0.98777,0.77847){\circle*{0.075}}
		\put(0.56163,0.12493){\circle*{0.075}}
		\put(0.55697,0.51967){\circle*{0.075}}
		\put(0.61643,0.84773){\circle*{0.075}}
		\put(0.15817,0.41287){\circle*{0.075}}
		\put(0.67923,0.85853){\circle*{0.075}}
		\put(0.51137,0.37807){\circle*{0.075}}
		\put(0.43003,0.63733){\circle*{0.075}}
		\put(0.53657,0.53527){\circle*{0.075}}
		\put(0.34883,0.46413){\circle*{0.075}}
		\put(0.35377,0.20447){\circle*{0.075}}
		\put(0.71563,0.41893){\circle*{0.075}}
		\put(0.28297,0.90567){\circle*{0.075}}
		\put(0.61043,0.38173){\circle*{0.075}}
		\put(0.84417,0.35887){\circle*{0.075}}
		\put(0.91323,0.03253){\circle*{0.075}}
		\put(0.75737,0.48407){\circle*{0.075}}
		\put(0.30403,0.85133){\circle*{0.075}}
		\put(0.94257,0.40127){\circle*{0.075}}
		\put(0.26283,0.11813){\circle*{0.075}}
		\put(0.51977,0.43047){\circle*{0.075}}
		\put(0.06963,0.91293){\circle*{0.075}}
		\put(0.80897,0.09167){\circle*{0.075}}
		\put(0.80443,0.11573){\circle*{0.075}}
		\put(0.33017,0.10487){\circle*{0.075}}
		\put(0.34723,0.40653){\circle*{0.075}}
		\put(0.80337,0.39007){\circle*{0.075}}
		\put(0.37803,0.26533){\circle*{0.075}}
		\put(0.14857,0.06727){\circle*{0.075}}
		\put(0.37683,0.97213){\circle*{0.075}}
		\put(0.48577,0.45647){\circle*{0.075}}
		\put(0.62363,0.60693){\circle*{0.075}}
		\put(0.13497,0.07767){\circle*{0.075}}
		\put(0.19843,0.04973){\circle*{0.075}}
		\put(0.61617,0.65087){\circle*{0.075}}
		\put(0.98123,0.98053){\circle*{0.075}}
		\put(0.64937,0.09607){\circle*{0.075}}
		\put(0.65203,0.87933){\circle*{0.075}}
		\put(0.15457,0.53327){\circle*{0.075}}
		\put(0.69083,0.02613){\circle*{0.075}}
		\put(0.25177,0.28247){\circle*{0.075}}
		\put(0.37763,0.50093){\circle*{0.075}}
		\end{picture}
	\begin{picture}(1.1,1.15)
\thicklines
\put(0,0){\line(1,0){1}}
\put(0,0){\line(0,1){1}}
\put(0,1){\line(1,0){1}}
\put(1,0){\line(0,1){1}}
\put(0.00097,0.00367){\circle*{0.15}}
\put(0.85243,0.84373){\circle*{0.15}}
\put(0.84217,0.53687){\circle*{0.15}}
\put(0.47523,0.01453){\circle*{0.15}}
\put(0.83537,0.54207){\circle*{0.15}}
\put(0.38603,0.55333){\circle*{0.15}}
\put(0.90057,0.13927){\circle*{0.15}}
\put(0.06483,0.74013){\circle*{0.15}}
\put(0.15777,0.64847){\circle*{0.15}}
\put(0.79163,0.65493){\circle*{0.15}}
\put(0.92697,0.58967){\circle*{0.15}}
\put(0.64643,0.17773){\circle*{0.15}}
\put(0.72817,0.68287){\circle*{0.15}}
\put(0.50923,0.98853){\circle*{0.15}}
\put(0.28137,0.84807){\circle*{0.15}}
\put(0.06003,0.56733){\circle*{0.15}}
\put(0.50657,0.20527){\circle*{0.15}}
\put(0.77883,0.19413){\circle*{0.15}}
\put(0.52377,0.07447){\circle*{0.15}}
\put(0.94563,0.94893){\circle*{0.15}}
\put(0.65297,0.97567){\circle*{0.15}}
\put(0.64043,0.71173){\circle*{0.15}}
\put(0.41417,0.62887){\circle*{0.15}}
\put(0.74323,0.16253){\circle*{0.15}}
\put(0.52737,0.95407){\circle*{0.15}}
\put(0.93403,0.78133){\circle*{0.15}}
\put(0.91257,0.07127){\circle*{0.15}}
\put(0.69283,0.84813){\circle*{0.15}}
\put(0.68977,0.30047){\circle*{0.15}}
\put(0.29963,0.44293){\circle*{0.15}}
\put(0.17897,0.16167){\circle*{0.15}}
\put(0.83443,0.44573){\circle*{0.15}}
\put(0.90017,0.37487){\circle*{0.15}}
\put(0.17723,0.53653){\circle*{0.15}}
\put(0.57337,0.86007){\circle*{0.15}}
\put(0.00803,0.19533){\circle*{0.15}}
\put(0.11857,0.73727){\circle*{0.15}}
\put(0.80683,0.70213){\circle*{0.15}}
\put(0.65577,0.32647){\circle*{0.15}}
\put(0.85363,0.13693){\circle*{0.15}}
\put(0.50497,0.14767){\circle*{0.15}}
\put(0.22843,0.37973){\circle*{0.15}}
\put(0.18617,0.92087){\circle*{0.15}}
\put(0.81123,0.11053){\circle*{0.15}}
\put(0.41937,0.56607){\circle*{0.15}}
\put(0.28203,0.80933){\circle*{0.15}}
\put(0.12457,0.20327){\circle*{0.15}}
\put(0.12083,0.75613){\circle*{0.15}}
\put(0.42177,0.15247){\circle*{0.15}}
\put(0.60763,0.03093){\circle*{0.15}}
\put(0.63097,0.93367){\circle*{0.15}}
\put(0.82243,0.51373){\circle*{0.15}}
\put(0.27217,0.26687){\circle*{0.15}}
\put(0.64523,0.88453){\circle*{0.15}}
\put(0.06537,0.07207){\circle*{0.15}}
\put(0.75603,0.62333){\circle*{0.15}}
\put(0.93057,0.46927){\circle*{0.15}}
\put(0.63483,0.01013){\circle*{0.15}}
\put(0.98777,0.77847){\circle*{0.15}}
\put(0.56163,0.12493){\circle*{0.15}}
\put(0.55697,0.51967){\circle*{0.15}}
\put(0.61643,0.84773){\circle*{0.15}}
\put(0.15817,0.41287){\circle*{0.15}}
\put(0.67923,0.85853){\circle*{0.15}}
\put(0.51137,0.37807){\circle*{0.15}}
\put(0.43003,0.63733){\circle*{0.15}}
\put(0.53657,0.53527){\circle*{0.15}}
\put(0.34883,0.46413){\circle*{0.15}}
\put(0.35377,0.20447){\circle*{0.15}}
\put(0.71563,0.41893){\circle*{0.15}}
\put(0.28297,0.90567){\circle*{0.15}}
\put(0.61043,0.38173){\circle*{0.15}}
\put(0.84417,0.35887){\circle*{0.15}}
\put(0.91323,0.03253){\circle*{0.15}}
\put(0.75737,0.48407){\circle*{0.15}}
\put(0.30403,0.85133){\circle*{0.15}}
\put(0.94257,0.40127){\circle*{0.15}}
\put(0.26283,0.11813){\circle*{0.15}}
\put(0.51977,0.43047){\circle*{0.15}}
\put(0.06963,0.91293){\circle*{0.15}}
\put(0.80897,0.09167){\circle*{0.15}}
\put(0.80443,0.11573){\circle*{0.15}}
\put(0.33017,0.10487){\circle*{0.15}}
\put(0.34723,0.40653){\circle*{0.15}}
\put(0.80337,0.39007){\circle*{0.15}}
\put(0.37803,0.26533){\circle*{0.15}}
\put(0.14857,0.06727){\circle*{0.15}}
\put(0.37683,0.97213){\circle*{0.15}}
\put(0.48577,0.45647){\circle*{0.15}}
\put(0.62363,0.60693){\circle*{0.15}}
\put(0.13497,0.07767){\circle*{0.15}}
\put(0.19843,0.04973){\circle*{0.15}}
\put(0.61617,0.65087){\circle*{0.15}}
\put(0.98123,0.98053){\circle*{0.15}}
\put(0.64937,0.09607){\circle*{0.15}}
\put(0.65203,0.87933){\circle*{0.15}}
\put(0.15457,0.53327){\circle*{0.15}}
\put(0.69083,0.02613){\circle*{0.15}}
\put(0.25177,0.28247){\circle*{0.15}}
\put(0.37763,0.50093){\circle*{0.15}}
\end{picture}
	\begin{picture}(1.1,1.15)
\thicklines
\put(0,0){\line(1,0){1}}
\put(0,0){\line(0,1){1}}
\put(0,1){\line(1,0){1}}
\put(1,0){\line(0,1){1}}
\put(0.00097,0.00367){\circle*{0.5}}
\put(0.85243,0.84373){\circle*{0.5}}
\put(0.84217,0.53687){\circle*{0.5}}
\put(0.47523,0.01453){\circle*{0.5}}
\put(0.83537,0.54207){\circle*{0.5}}
\put(0.38603,0.55333){\circle*{0.5}}
\put(0.90057,0.13927){\circle*{0.5}}
\put(0.06483,0.74013){\circle*{0.5}}
\put(0.15777,0.64847){\circle*{0.5}}
\put(0.79163,0.65493){\circle*{0.5}}
\put(0.92697,0.58967){\circle*{0.5}}
\put(0.64643,0.17773){\circle*{0.5}}
\put(0.72817,0.68287){\circle*{0.5}}
\put(0.50923,0.98853){\circle*{0.5}}
\put(0.28137,0.84807){\circle*{0.5}}
\put(0.06003,0.56733){\circle*{0.5}}
\put(0.50657,0.20527){\circle*{0.5}}
\put(0.77883,0.19413){\circle*{0.5}}
\put(0.52377,0.07447){\circle*{0.5}}
\put(0.94563,0.94893){\circle*{0.5}}
\put(0.65297,0.97567){\circle*{0.5}}
\put(0.64043,0.71173){\circle*{0.5}}
\put(0.41417,0.62887){\circle*{0.5}}
\put(0.74323,0.16253){\circle*{0.5}}
\put(0.52737,0.95407){\circle*{0.5}}
\put(0.93403,0.78133){\circle*{0.5}}
\put(0.91257,0.07127){\circle*{0.5}}
\put(0.69283,0.84813){\circle*{0.5}}
\put(0.68977,0.30047){\circle*{0.5}}
\put(0.29963,0.44293){\circle*{0.5}}
\put(0.17897,0.16167){\circle*{0.5}}
\put(0.83443,0.44573){\circle*{0.5}}
\put(0.90017,0.37487){\circle*{0.5}}
\put(0.17723,0.53653){\circle*{0.5}}
\put(0.57337,0.86007){\circle*{0.5}}
\put(0.00803,0.19533){\circle*{0.5}}
\put(0.11857,0.73727){\circle*{0.5}}
\put(0.80683,0.70213){\circle*{0.5}}
\put(0.65577,0.32647){\circle*{0.5}}
\put(0.85363,0.13693){\circle*{0.5}}
\put(0.50497,0.14767){\circle*{0.5}}
\put(0.22843,0.37973){\circle*{0.5}}
\put(0.18617,0.92087){\circle*{0.5}}
\put(0.81123,0.11053){\circle*{0.5}}
\put(0.41937,0.56607){\circle*{0.5}}
\put(0.28203,0.80933){\circle*{0.5}}
\put(0.12457,0.20327){\circle*{0.5}}
\put(0.12083,0.75613){\circle*{0.5}}
\put(0.42177,0.15247){\circle*{0.5}}
\put(0.60763,0.03093){\circle*{0.5}}
\put(0.63097,0.93367){\circle*{0.5}}
\put(0.82243,0.51373){\circle*{0.5}}
\put(0.27217,0.26687){\circle*{0.5}}
\put(0.64523,0.88453){\circle*{0.5}}
\put(0.06537,0.07207){\circle*{0.5}}
\put(0.75603,0.62333){\circle*{0.5}}
\put(0.93057,0.46927){\circle*{0.5}}
\put(0.63483,0.01013){\circle*{0.5}}
\put(0.98777,0.77847){\circle*{0.5}}
\put(0.56163,0.12493){\circle*{0.5}}
\put(0.55697,0.51967){\circle*{0.5}}
\put(0.61643,0.84773){\circle*{0.5}}
\put(0.15817,0.41287){\circle*{0.5}}
\put(0.67923,0.85853){\circle*{0.5}}
\put(0.51137,0.37807){\circle*{0.5}}
\put(0.43003,0.63733){\circle*{0.5}}
\put(0.53657,0.53527){\circle*{0.5}}
\put(0.34883,0.46413){\circle*{0.5}}
\put(0.35377,0.20447){\circle*{0.5}}
\put(0.71563,0.41893){\circle*{0.5}}
\put(0.28297,0.90567){\circle*{0.5}}
\put(0.61043,0.38173){\circle*{0.5}}
\put(0.84417,0.35887){\circle*{0.5}}
\put(0.91323,0.03253){\circle*{0.5}}
\put(0.75737,0.48407){\circle*{0.5}}
\put(0.30403,0.85133){\circle*{0.5}}
\put(0.94257,0.40127){\circle*{0.5}}
\put(0.26283,0.11813){\circle*{0.5}}
\put(0.51977,0.43047){\circle*{0.5}}
\put(0.06963,0.91293){\circle*{0.5}}
\put(0.80897,0.09167){\circle*{0.5}}
\put(0.80443,0.11573){\circle*{0.5}}
\put(0.33017,0.10487){\circle*{0.5}}
\put(0.34723,0.40653){\circle*{0.5}}
\put(0.80337,0.39007){\circle*{0.5}}
\put(0.37803,0.26533){\circle*{0.5}}
\put(0.14857,0.06727){\circle*{0.5}}
\put(0.37683,0.97213){\circle*{0.5}}
\put(0.48577,0.45647){\circle*{0.5}}
\put(0.62363,0.60693){\circle*{0.5}}
\put(0.13497,0.07767){\circle*{0.5}}
\put(0.19843,0.04973){\circle*{0.5}}
\put(0.61617,0.65087){\circle*{0.5}}
\put(0.98123,0.98053){\circle*{0.5}}
\put(0.64937,0.09607){\circle*{0.5}}
\put(0.65203,0.87933){\circle*{0.5}}
\put(0.15457,0.53327){\circle*{0.5}}
\put(0.69083,0.02613){\circle*{0.5}}
\put(0.25177,0.28247){\circle*{0.5}}
\put(0.37763,0.50093){\circle*{0.5}}
\end{picture}
		\caption{\label{points} Configuration of 100 communities in $N_o$ and various values of $D$}
	\end{center}
\end{figure}
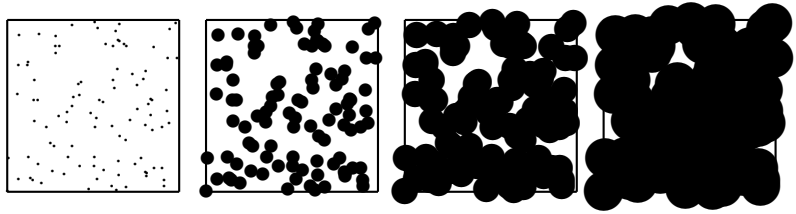

\begin{figure}[ht!]
	\begin{center}
		\setlength{\unitlength}{1in}
		\begin{picture}(5,5.1)
		\includegraphics[width=5in]{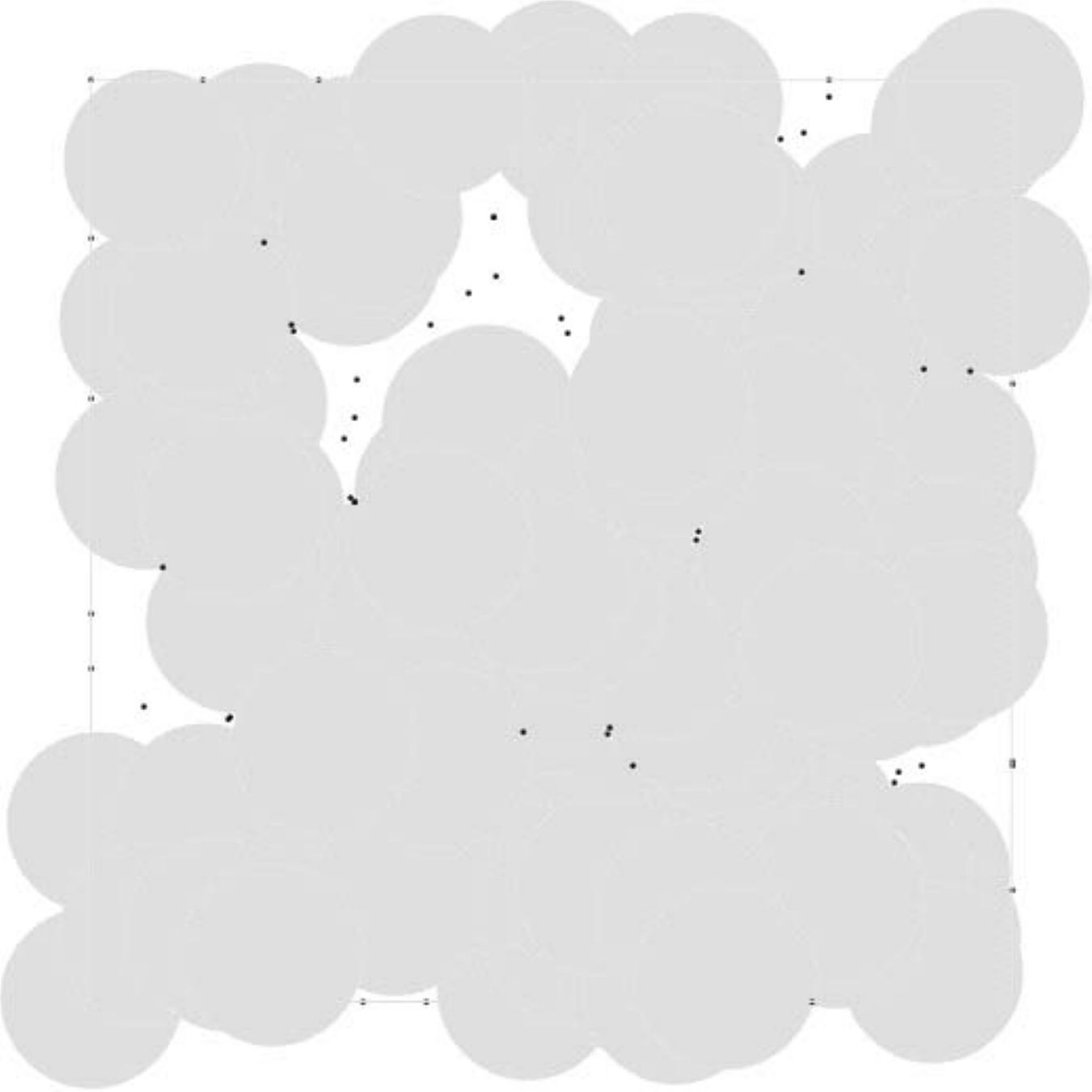}
		\end{picture}
		\caption{\label{50}Infeasible areas for $n=100$, $D=0.95$, and 50 Voronoi points}
	\end{center}
\end{figure}

In order to simplify the presentation we assume in the following discussion, unless specified otherwise,  that $N_o=N_d=N$ meaning that every point generates demand and is also negatively affected by the facilities. In Figure \ref{points} increasing values of $D$ are depicted.
When $D=0$, no distance constraints are included in the model and the obnoxious $p$-median problem is the standard one. As $D$ increases, small circles surrounding each demand point form infeasible regions. If $D$ is quite small the feasible region is well connected and solving the unconstrained $p$-median problem may constitute a good starting solution for a follow-up phase incorporating the obnoxious constraints. When $D$ is large, the whole convex hull of the demand points becomes infeasible and the solution is to locate the facilities at the periphery of the infeasible region.

The most challenging problem seems to be moderate values of $D$. As $D$ increases, the feasible region decreases and becomes a union of disconnected areas. If a starting solution includes a facility in one of these disconnected feasible areas, standard non-linear procedures will not move such a facility to another feasible area. To get a ``good" solution the feasible areas for the facilities need to be properly selected.

\section{\label{sec3}Finding the Feasible Areas}

\begin{figure}[ht!]
	\begin{center}
		\setlength{\unitlength}{1in}
		\begin{picture}(5,5.1)
		\includegraphics[width=5in]{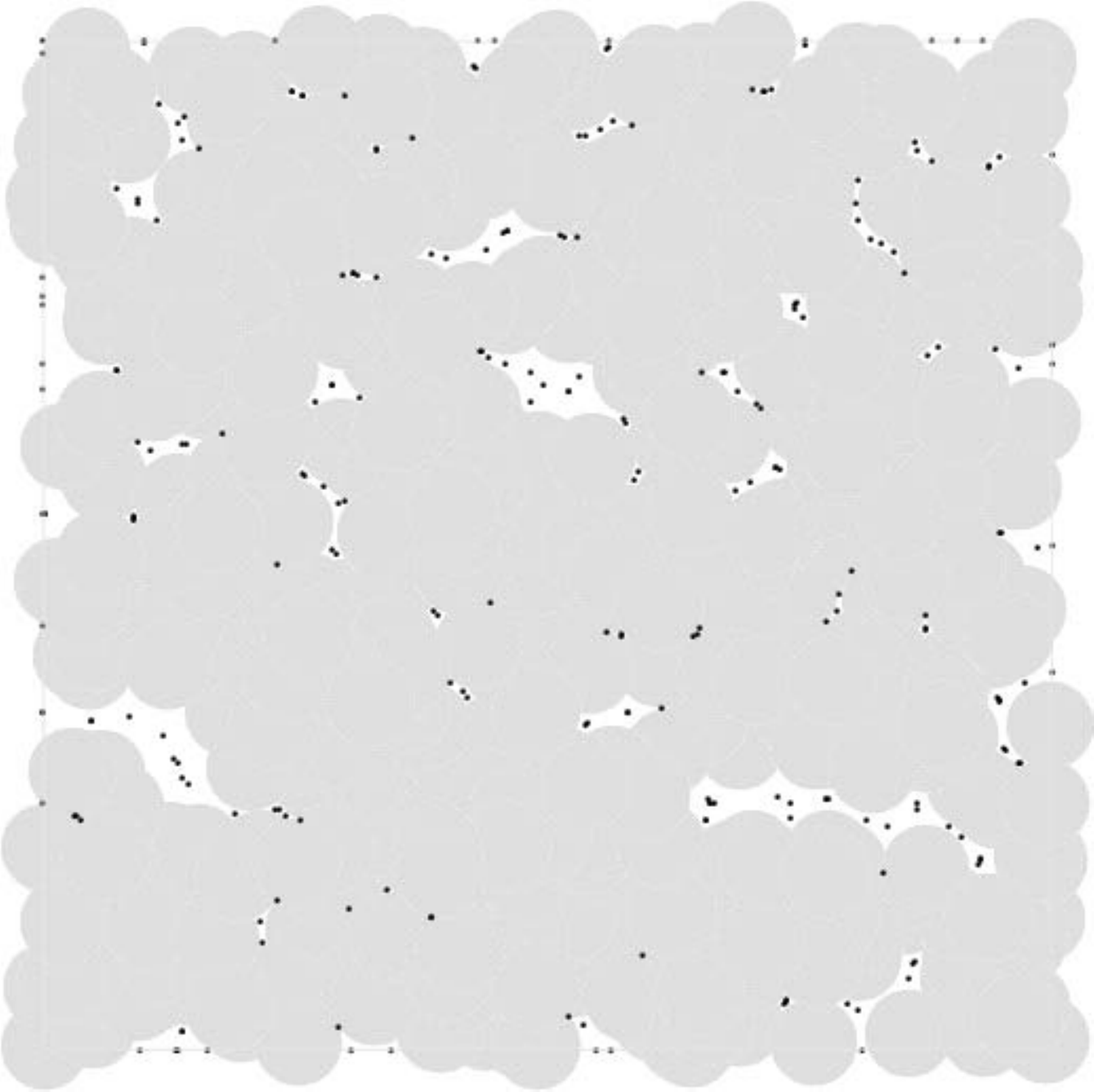}
		\end{picture}
		\caption{\label{500}Infeasible areas for $n=500$, $D=0.42$, and 245 Voronoi points}
	\end{center}
\end{figure}

\begin{figure}[ht!]
	\begin{center}
		\setlength{\unitlength}{1in}
		\begin{picture}(5,5.1)
		\includegraphics[width=5in]{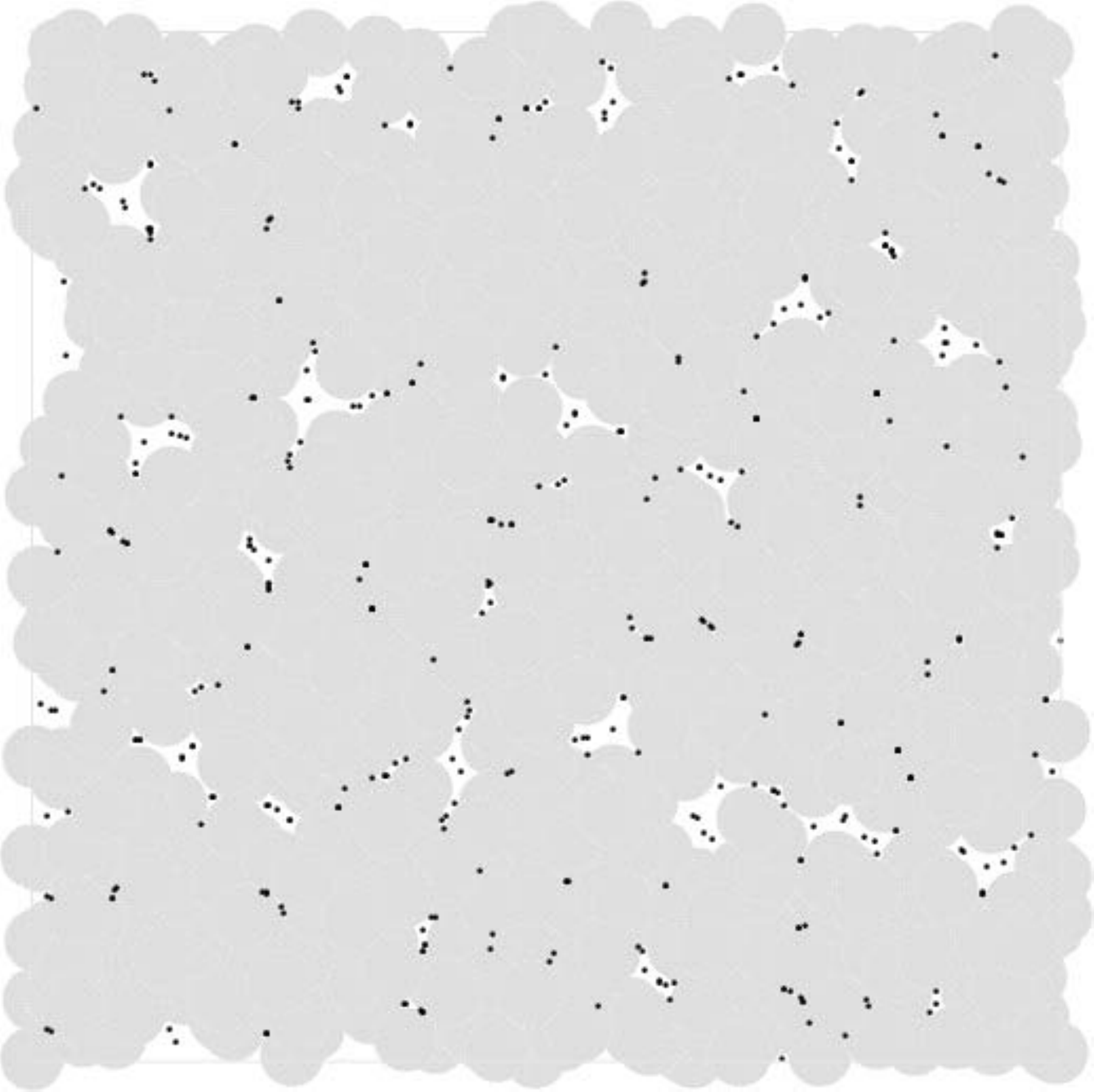}
		\end{picture}
				\caption{\label{1000}Infeasible areas for $n=1000$, $D=0.3$, and 473 Voronoi points}
			\end{center}
\end{figure}

\begin{table}[ht!]
	\begin{center}
		\caption{\label{Vor}The first 50 Voronoi points ($n=100$)}
		\medskip
		
		\begin{tabular}{|c|c|c|c||c|c|c|c|}
			\hline
			$i$&$x_i$&$y_i$&$D_i$&$i$&$x_i$&$y_i$&$D_i$\\
			\hline
			1&	0.00000&	3.61453&	1.66317&	26&	0.00000&	6.55464&	1.06636\\
			2&	0.00000&	4.20781&	1.58368&	27&	1.48262&	3.06965&	1.06367\\
			3&	10.00000&	2.57239&	1.54282&	28&	2.81536&	5.46099&	1.04744\\
			4&	8.02745&	10.00000&	1.51738&	29&	7.49008&	9.35686&	1.04029\\
			5&	4.40903&	7.87825&	1.50887&	30&	8.72091&	2.38579&	1.03312\\
			6&	10.00000&	2.61785&	1.50845&	31&	1.52648&	3.10114&	1.02905\\
			7&	8.01192&	9.83008&	1.48404&	32&	6.58788&	5.10980&	1.02189\\
			8&	0.57979&	3.21438&	1.35640&	33&	6.57244&	5.02285&	1.01770\\
			9&	4.38806&	8.52444&	1.34780&	34&	1.22238&	10.00000&	1.01729\\
			10&	4.38830&	8.52487&	1.34754&	35&	2.19860&	7.28931&	1.01632\\
			11&	2.88799&	6.75677&	1.33824&	36&	5.61395&	2.90214&	1.01100\\
			12&	4.09488&	7.69328&	1.33587&	37&	10.00000&	1.21758&	1.00960\\
			13&	5.11154&	7.42616&	1.32914&	38&	2.85770&	5.42785&	1.00813\\
			14&	2.85594&	6.33755&	1.28668&	39&	2.48189&	10.00000&	1.00538\\
			15&	7.75331&	9.43913&	1.26415&	40&	2.86769&	5.42369&	0.99864\\
			16&	2.76281&	6.11411&	1.24170&	41&	9.04853&	6.86445&	0.99270\\
			17&	5.17172&	7.25598&	1.24036&	42&	9.55826&	6.84567&	0.99187\\
			18&	9.03540&	2.57075&	1.17843&	43&	2.17176&	7.35019&	0.98631\\
			19&	3.68326&	7.33910&	1.14609&	44&	0.78639&	4.70746&	0.98361\\
			20&	8.77401&	2.50362&	1.13482&	45&	5.89198&	2.55987&	0.96952\\
			21&	0.00000&	10.00000&	1.11488&	46&	7.71137&	7.91767&	0.96482\\
			22&	3.65520&	0.00000&	1.10668&	47&	5.63950&	2.97363&	0.96324\\
			23&	2.96654&	0.00000&	1.10096&	48&	1.87471&	8.25026&	0.95853\\
			24&	0.00000&	8.28398&	1.09517&	49&	7.82576&	0.00000&	0.95394\\
			25&	10.00000&	6.70342&	1.08818&	50&	4.69156&	2.92776&	0.95169\\
			\hline							
		\end{tabular}
	\end{center}
\end{table}

\citet{DDS18a} proposed the ``Big Arc Small Arc" algorithm for solving the single facility location problem. All feasible arcs on the periphery of circles are found and evaluated by a global optimization algorithm. A lower bound is found on each feasible arc and if the best found solution is higher than the lower bound, the arc is divided into two arcs and the iterations continue until the difference between the best found solution and the lowest lower bound is less than $\epsilon>0$. This approach is similar to many global optimization algorithms widely used to solve non-convex location problems: Big Square Small Square \citep{HPT81}, Big Region Small Region \citep{HJK95}, Big Triangle Small Triangle \citep{DS04,S19}, Big Cube Small Cube \citep{SS10}, Big Segment Small Segment \citep{BDK09seg}.
If the optimal location of the facility is not feasible (if it is feasible it is optimal for the constrained problem), then the optimal solution must be on a feasible arc of a circle. However, for multiple facilities this property does not necessarily hold because the problem is not convex. Some facilities may be located on feasible arcs and some may be located in ``open" feasible areas. We therefore propose a different approach for heuristically solving the multiple facility version of this problem.
 
In order to find the distinct feasible areas we employ the concept of the Voronoi diagram \citep{SOK95,OKSC00,V08,AKL13}  based on the points in $N_o$. The plane is partitioned into triangles \citep[termed Delaunay triangulation; ][]{LS80} so that the centers of the triangles are equidistant from the three vertices of the triangles and are termed Voronoi points. These Voronoi points can be found by Mathematica \citep{Math} and other available software \citep{OIM84,SI92}. Each point (a vertex of a triangle) is closest to a Voronoi point located at that triangle. There are $V$ Voronoi points.  Voronoi point $V_j$ is at distance $D_j$ from the closest demand point in $N_o$ for $j=1,\ldots V$. The Voronoi points can be sorted from the largest minimum distance down creating a list of distances $D_1\ge D_2\ge\ldots\ge D_V$. \citet{Sham} used this scheme to find the location for a facility which is farthest from all demand points by choosing the Voronoi point with minimum distance $D_1$. \citet{DKS18} used the list of Voronoi points to heuristically solve a multiple obnoxious facilities problem.

 In Figures \ref{50}-\ref{1000} the infeasible areas are depicted for the problems used in the computational experiments with $n=100$ demand points and $D=0.95$, $n=500$ and $D=0.42$, and $n=1000$ with $D=0.3$. The Voronoi points are marked by black dots. One may have the impression that many Voronoi points are inside the gray  infeasible regions. However, they are located in very small feasible regions so that the white area is smaller than the black dot. It will be very difficult for conventional approaches to identify such small feasible areas without identifying the Voronoi points.

\begin{figure}[ht!]
	\begin{center}
		\setlength{\unitlength}{1in}
		\begin{picture}(2,2.)
		\includegraphics[width=2in]{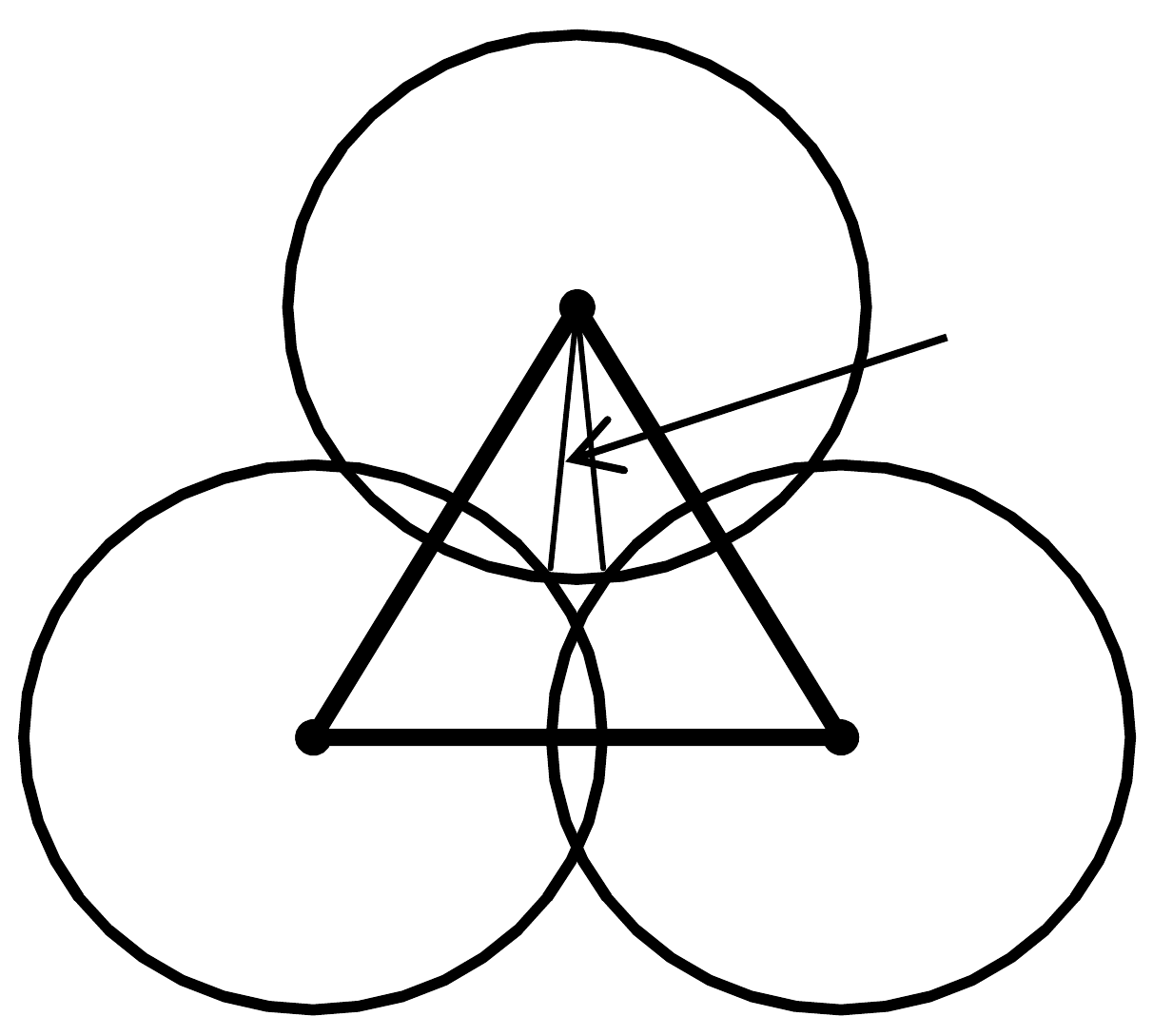}
		\put(-0.35,1.16){$2\theta$}
		\end{picture}
		\caption{\label{triang}The triangle example}
	\end{center}
\end{figure}

To illustrate this phenomenon consider an equilateral triangle, whose vertices are demand points, depicted in Figure \ref{triang} for which the distance to the Voronoi point $D_i$ is $1.05D$.
The radii of the circles (infeasible areas) are $D$. The feasible area is at the center of the triangle outside the circles. The bottom side of the long narrow triangle is $2D\sin\theta$. The distance to the center of the tiny triangle (the Voronoi point) is at distance $D_i$ from the top vertex. Therefore,
$D(1-\cos\theta)+D_i-D=\frac{1}{\sqrt3}D\sin\theta$ leading to the relationship $\frac{1}{\sqrt3}\sin\theta +\cos\theta=\frac{D_i}{D}$. Multiplying by $\frac{\sqrt3}{2}$ and solving for $\theta$: 
$$\frac12\sin\theta+\frac{\sqrt{3}}2\cos\theta=\frac{\sqrt3D_i}{2D}~\rightarrow~
\sin(\frac\pi3+\theta)=\frac{\sqrt3D_i}{2D}~\rightarrow~\theta=\arcsin\left\{\frac{\sqrt3 D_i}{2D}\right\}-\frac\pi3~.$$
Note that if $D_i>\frac2{\sqrt3} D$ the circles do not intersect.

The area of the small triangle is $\sqrt3 D^2\sin^2\theta$.
The area of the part of the circle which is inside the triangle (non-feasible) is $D^2\theta-D^2\sin\theta\cos\theta$.
The area $A$ of the feasible region is therefore
\begin{eqnarray}
A&=&\sqrt3 D^2\sin^2\theta-3[D^2\theta-D^2\sin\theta\cos\theta]=
D^2\left\{\sqrt3 \sin^2\theta-3[\theta-\sin\theta\cos\theta]\right\}\nonumber\\
&=&2\sqrt3D^2\sin\theta[\frac12\sin\theta+\frac{\sqrt3}{2}\cos\theta]-3\theta D^2=3D[D_i\sin\theta-D\theta]~.
\label{arexact}
\end{eqnarray}
Assuming a small $\theta$,  $\theta\approx\sin\theta \approx \sqrt3\frac{D_i-D}{D}$ and thus by approximating equation (\ref{arexact}):

\begin{equation}\label{ar}
A\approx 3\sqrt3 (D_i-D)^2~.
\end{equation} 

For example, for $D=1$ $D_i=1.05$, the exact area is 0.013727 and the approximated area is 0.012990. For $D=1$ $D_i=1.01$ the areas are 0.000525 and 0.000520. The exact results were confirmed by simulating randomly generated one billion points.

The largest distance $d_{\max}$ between the Voronoi point and the feasible area  is
\begin{equation}\label{dis}
d_{\max}=\frac{D}{\sqrt3}\sin\theta\approx D_i-D~. 
\end{equation}
For $D=0.95$ and $n=100$ there are 10 Voronoi points satisfying $0.95<D_i< 0.95\times 1.05=0.9975$ (see Table \ref{Vor}). These points are at most at distance 0.05 from the boundary of their feasible area and clearly appear as black points surrounded by a gray area in Figure~\ref{50}. Recall that the square size is 10 by 10,  and such a distance is $\frac{1}{200}$ of the square's side. The side of the square in the figure is about 4 inches and thus this distance is 0.02 inches or 0.5 millimeters.

\section{\label{sec4}The Proposed Solution Algorithm}

We propose to find all the Voronoi points, identify which ones are feasible, find the $p$ feasible Voronoi points that yield the best value of the objective function, and use these points as a starting solution for non-linear solvers. 

All the Voronoi points for which $D_i\ge D$ form a region of feasible area (if $D_i=D$ the region is one point). Actually, if $D> D_1$, the whole convex hull of demand points is infeasible. Suppose that the first $m\ge p$ sorted Voronoi points satisfy $D_i\ge D$. We propose to find the $p$-median solution for facilities located at $p$ of these $m$ Voronoi points and use it as a starting (feasible) solution for non-linear optimization procedures.
This is a discrete $p$-median problem that can be solved by the following Mixed Binary Linear Program (MBLP) or Binary Linear Program (BLP) \citep{Das95}.

Let $d_{ij}$ is the distance between demand point $i$ and potential location (Voronoi point) $j$.
We minimize the value of the objective function by the best selection of the set $P$ of the $p$ out of the $m$ potential locations. The standard formulation for the discrete $p$-median problem is:
\begin{equation}
\min\limits_{P}\left\{~\sum\limits_{i\in N_d}w_i\min\limits_{j\in P}\{d_{ij}\}~\right\}
\end{equation}

Let $x_j\in \{0,1\}$ be a binary variable. $x_j=1$ if location $j$ is selected and zero otherwise. For MBLP $0\le y_{ij}\le 1$ are continuous variables for $i\in N_d$, $j=1,\ldots,m$. For BLP $y_{ij}\in\{0,1\}$. The MBLP/BLP optimization problem is:
\begin{eqnarray}
&\min\left\{\sum\limits_{i\in N_d}\sum\limits_{j=1}^m\left[w_id_{ij}\right]y_{ij}\right\}\label{BLP}\\
\mbox{subject to:}&&\nonumber\\
&\sum\limits_{j=1}^mx_j=p\nonumber\\
&y_{ij}\le x_j& \mbox{ for }i\in N_d;~j=1,\ldots m\nonumber\\
&\sum\limits_{j=1}^my_{ij}=1& \mbox{ for }i\in N_d\nonumber\\
&x_j\in\{0,1\}\nonumber\\
\mbox{For MBLP:}&y_{ij}\ge 0\nonumber\\
\mbox{For BLP:}&y_{ij}\in\{0,1\}\nonumber
\end{eqnarray}
Note that this formulation can be used for any distance measure. In the solution $y_{ij}\in\{0,1\}$ because for each demand point $y_{ij}=1$ for the closest facility and all others are zeros. If there is a tie in the minimum distance, $y_{ij}$ can be non-integer but the objective function is the same as the integer solution. Therefore, an equivalent formulation is to require that $y_{ij}$ are binary variables yielding a BLP formulation. The MBLP and BLP formulations were solved by CPLEX \citep{cplex}. The BLP performed faster and thus was used in the computational experiments.

We concentrate on medium values of $D$ which are the most difficult to solve. We therefore ignore the possibility that a facility will be located outside the convex hull. It may be an issue if $m$ is not much larger than $p$. It is of course an issue if $m<p$. It is possible to add all feasible intersection points which are outside the convex hull of demand points. However, this may increase the size of the BLP (\ref{BLP}) so it will have to be solved heuristically rather than optimally.

\section{\label{comp}Computational Experiments}

In all the computational experiments it is assumed that $N_d=N_o=N$. The number of demand points is denoted by $n$. CPLEX \citep{cplex} and SNOPT \citep{SNOPT} were run on an Amazon EC2 instance with 32 CPUs and 70 GB of RAM. The BLP was implemented with the OPL and run on IBM's CPLEX Optimization Studio 12.8 environment. We used the default CPLEX MIP solver settings.
The unconstrained $p$-median problem algorithm was coded in Fortran using double precision arithmetic. The programs were compiled by an Intel 11.1 Fortran Compiler with no parallel processing. They were run on a desktop with the Intel i7-6700 3.4GHz CPU processor and 16GB RAM.

To allow for easy replication in future comparisons, we tested randomly generated instances by the method proposed in \citet{DKS18,DDK19}.
A sequence of integer numbers in the open range (0, 100,000) is generated. A starting seed $r_1$, which is the first number in the sequence, is selected. The sequence is generated by the following rule for $k\ge 1$:
\begin{itemize}\narrow
	\item Set $\theta=12219 r_k$.
	\item Set $r_{k+1}=\theta-\lfloor\frac{\theta}{100000}\rfloor\times 100000$, i.e., $r_{k+1}$ is the remainder of dividing $\theta$ by 100000.
\end{itemize}

For the $x$ coordinates we used $r_1=97$ and for the $y$-coordinates we used $r_1=367$.  To define the coordinates in a 10 by 10 square, $r_k$ is divided by 10000.

The 50 Voronoi points (out of the $V=202$ Voronoi points) with the largest values of $D_i$ for the $n=100$ problem are depicted in Table \ref{Vor}. These  are all the Voronoi points whose $D_i\ge 0.95$ miles (in a 10 by 10 square mile). Recall that all Voronoi points for which $D_i\ge D$ are feasible. 

For the p-median problem we use all $w_i=1$.  1000 points are generated and the first $n$ points selected. Three value of $n$ were tested, $n=100,500,1000$ with corresponding values of $D=0.95,0.42,0.3$. These values of $D$ yielded 50, 245, and 473 feasible Voronoi points, respectively. See Figures~\ref{50}-\ref{1000}.

We first solved the unconstrained $p$-median for $p=2,3,4,5,10,15,20$, and $n=100,500,1000$ demand points, by the best available heuristic solution method \citep{DTD18}. Each instance was solved 10 times from randomly generated starting solutions and the best solution was found in all 10 replications for every instance. Then we found by solving BLP (\ref{BLP}) the optimal obnoxious $p$ median solutions when facilities are restricted to the set of feasible Voronoi points. These solutions where used as starting solutions for applying SNOPT \citep{SNOPT} once. The SNOPT results using the BLP results as starting solutions  serve as bench mark values for the obnoxious $p$-median solutions.

We then applied SNOPT from 1000 randomly generated starting solutions for each instance. This process took about 1000 times longer and yielded very poor results.  These results of random starting solutions (feasible or not) are not reported. 
We therefore randomly generated 100,000 random points and selected a subset of these getting a set of feasible points for each $n$. We then solved the obnoxious $p$-median problem by SNOPT using randomly selected $p$ points from the subset of feasible points as starting solutions. This approach produced significantly better results than those obtained by randomly generating 1000 starting solutions, whether feasible or not. However, run times are still 1000 longer than starting from Voronoi points.

\begin{table}[ht!]
	\begin{center}
		\caption{\label{Res}Objective Function Results}
		\medskip
		
	\setlength{\tabcolsep}{8pt}
		\begin{tabular}{|c|c|c||c|c|c|}
			\hline
			$p$&(1)&(2)&(3)&(4)&(5)\\
			\hline
\multicolumn{6}{|c|}{$n=100$}\\
\hline
2&	288.99&	293.66&	292.62&	292.62&	0.00\%\\
3&	228.67&	242.10&	241.15&	242.20&	0.44\%\\
4&	192.25&	209.54&	207.52&	210.63&	1.50\%\\
5&	164.60&	188.00&	185.80&	191.06&	2.83\%\\
10&	100.76&	142.60&	139.40&	148.54&	6.56\%\\
15&	74.47&	131.57&	126.23&	133.79&	5.99\%\\
20&	59.48&	127.48&	119.48&	133.69&	11.89\%\\
			\hline
\multicolumn{6}{|c|}{$n=500$}\\
\hline
2&	1476.06&	1501.01&	1497.92&	1497.96&	0.00\%\\
3&	1154.12&	1175.26&	1169.50&	1181.80&	1.05\%\\
4&	950.36&	965.45&	964.11&	974.19&	1.05\%\\
5&	856.12&	879.95&	874.82&	908.34&	3.83\%\\
10&	575.67&	619.30&	614.91&	667.67&	8.58\%\\
15&	449.89&	515.68&	508.18&	574.97&	13.14\%\\
20&	382.69&	452.57&	445.04&	510.57&	14.72\%\\
\hline
\multicolumn{6}{|c|}{$n=1000$}\\
\hline
2&	2912.11&	2945.71&	2943.87&	2948.43&	0.15\%\\
3&	2304.30&	2324.58&	2323.52&	2367.86&	1.91\%\\
4&	1920.37&	1922.81&	1921.37&	1954.25&	1.71\%\\
5&	1731.63&	1752.09&	1750.51&	1798.90&	2.76\%\\
10&	1177.97&	1219.96&	1215.08&	1367.36&	12.53\%\\
15&	942.47&	993.98&	988.84&	1110.75&	12.33\%\\
20&	798.55&	868.66&	863.50&	996.90&	15.45\%\\
\hline
\multicolumn{6}{l}{(1) For reference: unconstrained solution.}\\
\multicolumn{6}{l}{(2) For reference: best value at Voronoi points.}\\
\multicolumn{6}{l}{(3) SNOPT result run once from (2) solution.}\\
\multicolumn{6}{l}{(4) Best of 1000 SNOPT runs (random feasible points).}\\
\multicolumn{6}{l}{(5) \%  of (4) solution above (3) solution.}\\
						
		\end{tabular}
	\end{center}
\end{table}

The results of these experiments are depicted in Table \ref{Res}. Using  the BLP (\ref{BLP}) results as starting solutions performed best both in terms of run times and quality of results. Run times (not reported) are about 1000 times shorter. The random results matched the BLP starting solution results in only one instance out of 21 instances and are as much as 15\% higher.
 
The reason for the under-performance of randomly selecting feasible points is that there are many very small feasible areas. The black dots  in Figures \ref{50}-\ref{1000} represent Voronoi points which are all feasible but many look like they are located in an infeasible gray area. When a Voronoi point distance to the closest demand point is $D_i$, see Figure \ref{triang}, the feasible area by equation (\ref{ar}) is approximately $5(D_i-D)^2$  . For $D=0.95$ and $D_i=0.96$ the area is about $5\times 10^{-4}$ which is about 1/200,000 of the area of the 10 by 10 square. There are three Voronoi points for which $D_i<0.96$, see Table~\ref{Vor}. The feasible area has an area of about 5. Therefore, the probability that the generated feasible point is in that small area is about $ 10^{-4}$. Even when 1,000 random feasible points are generated, the chance that even one point is in such a feasible area is very low (9.5\%). Therefore, if a good solution has some facilities located in such small areas, it is likely to be missed even when 1000 random feasible starting solutions are used. Even if a particular small area is selected, it is unlikely that the additional required  areas for locating other facilities will all be selected as well. When the Voronoi points are used, there is a candidate location in every feasible area and the BLP selects the best set of candidate locations. Note also that if the good solution has a point in a small area, the final location of this facility is very close to the Voronoi point, equation (\ref{dis}), and thus the value of the objective function changes very little and the ``correct" feasible regions are likely to be selected by the BLP.

\section{\label{Efficient}The Efficient Frontier}

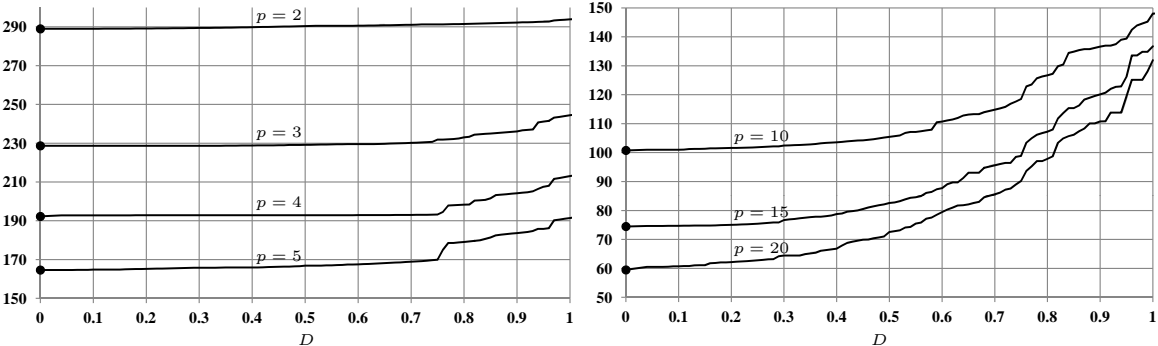
\begin{figure}[ht!]
	\begin{center}
		\setlength{\unitlength}{1in}
		\begin{picture}(6,1.9)	
		\includegraphics[width=3in]{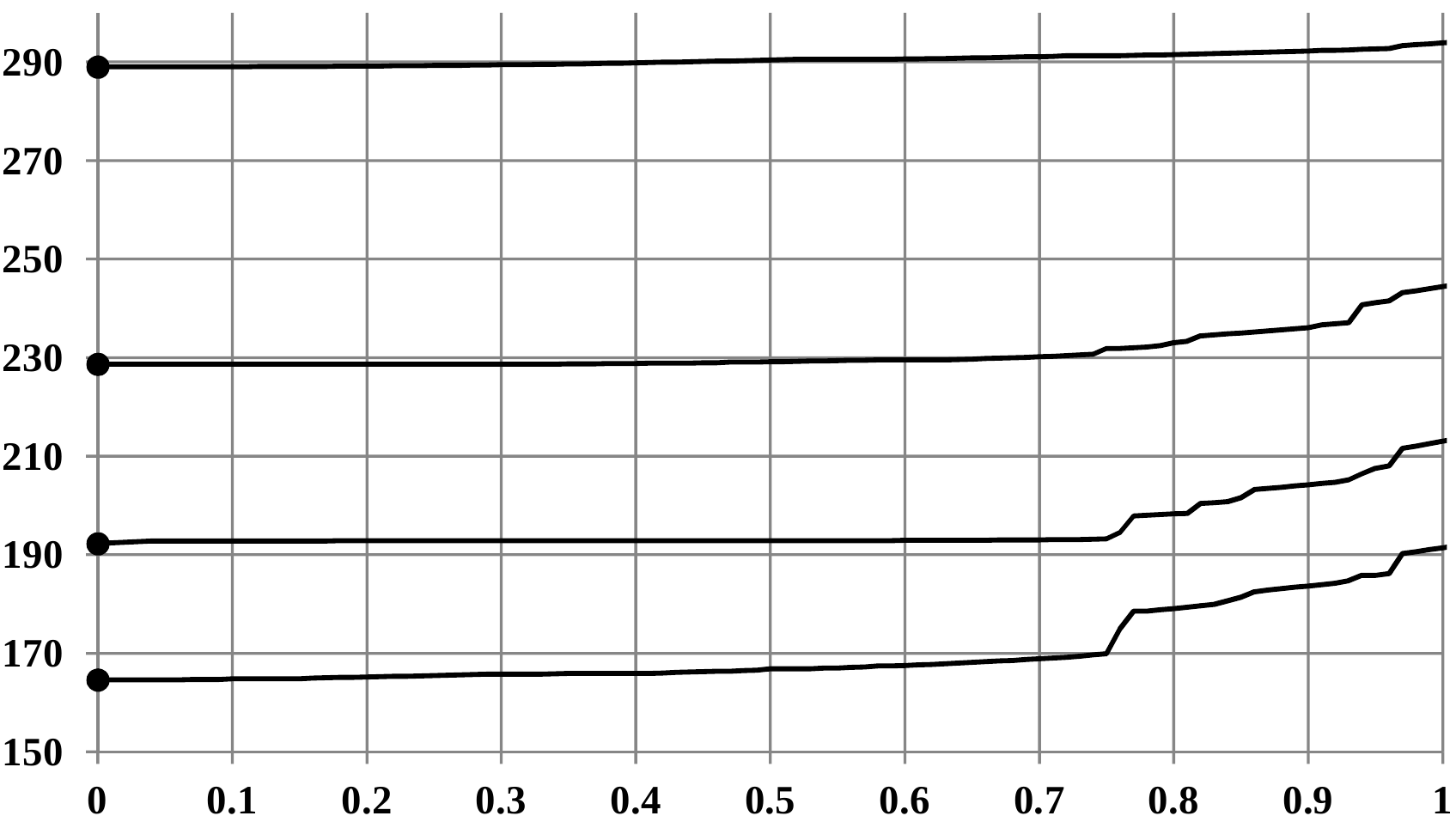}
		\includegraphics[width=3in]{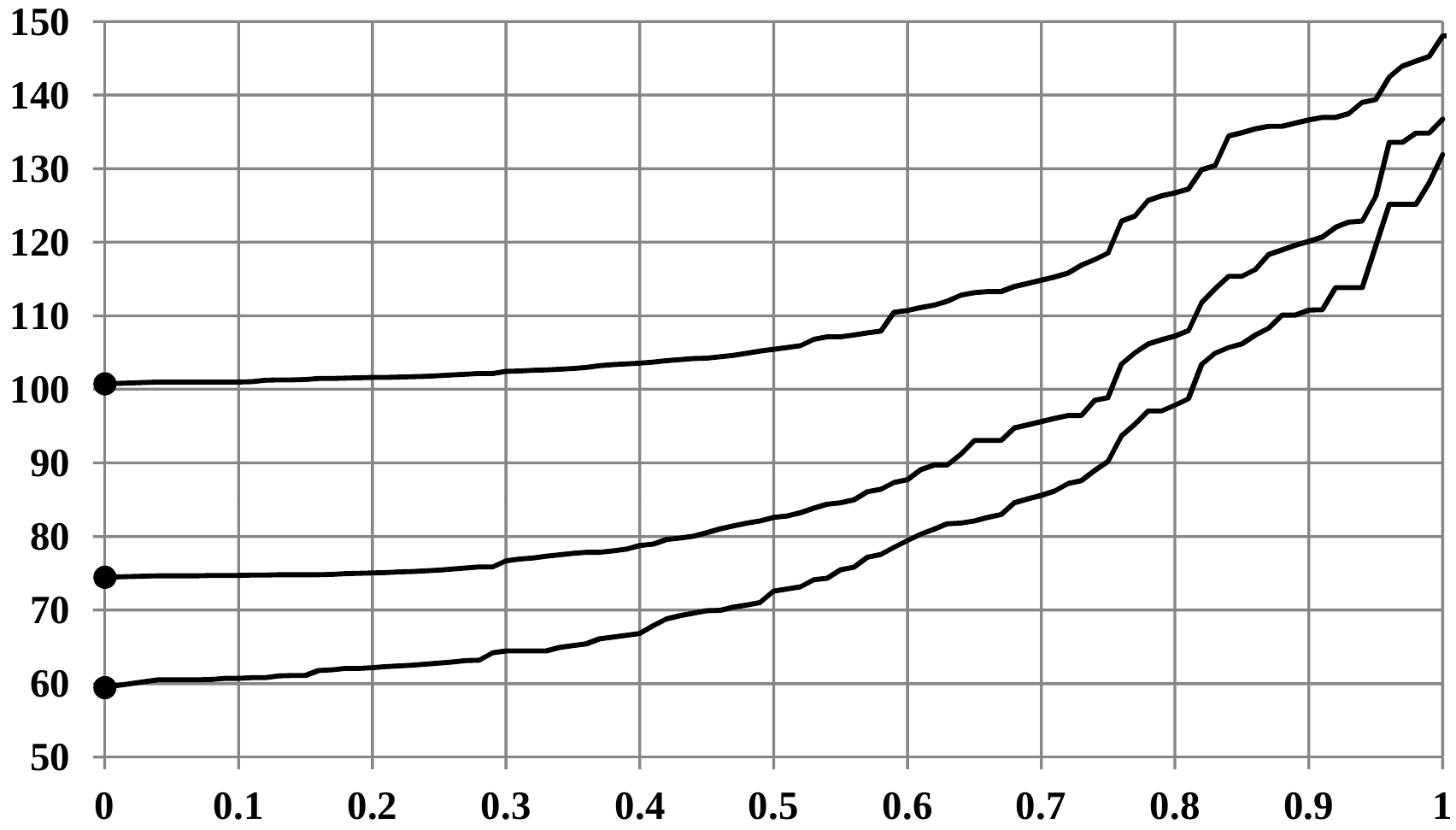}
		\scriptsize
		\tiny
		\put(-4.5,-0.1){$D$}
		\put(-1.5,-0.1){$D$}
		\put(-4.72,1.6){$p=2$}
		\put(-4.72,0.99){$p=3$}
		\put(-4.72,0.62){$p=4$}
		\put(-4.72,0.34){$p=5$}
		\put(-2.22,0.97){$p=10$}
		\put(-2.22,0.57){$p=15$}
		\put(-2.22,0.38){$p=20$}
		\end{picture}
		\caption{\label{ef100}The Efficient Frontiers for $n=100$}
	\end{center}
\end{figure}

\begin{figure}[ht!]
	\begin{center}
		\setlength{\unitlength}{1in}
		\begin{picture}(6,1.9)	
		\includegraphics[width=3in]{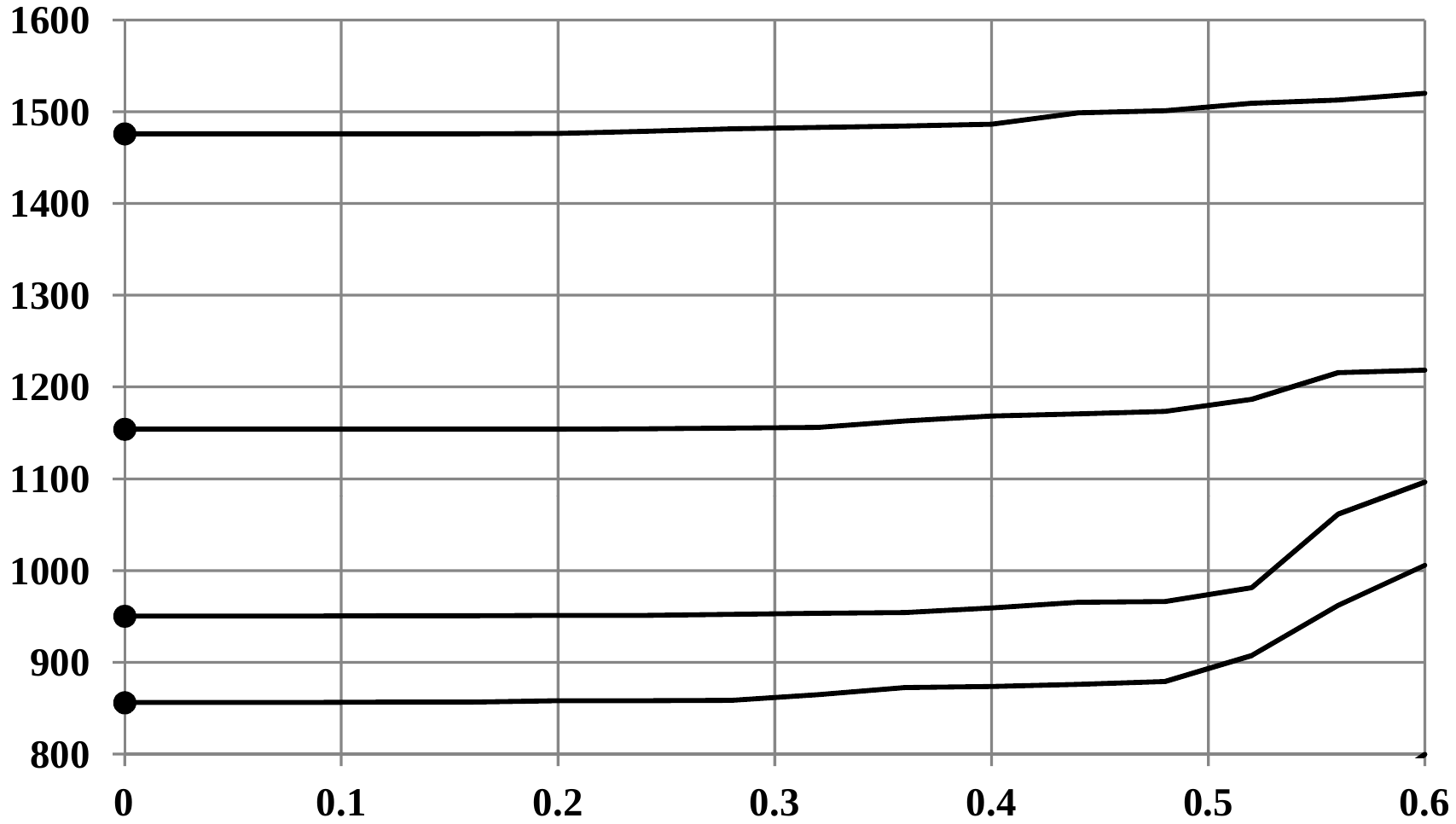}
		\includegraphics[width=3in]{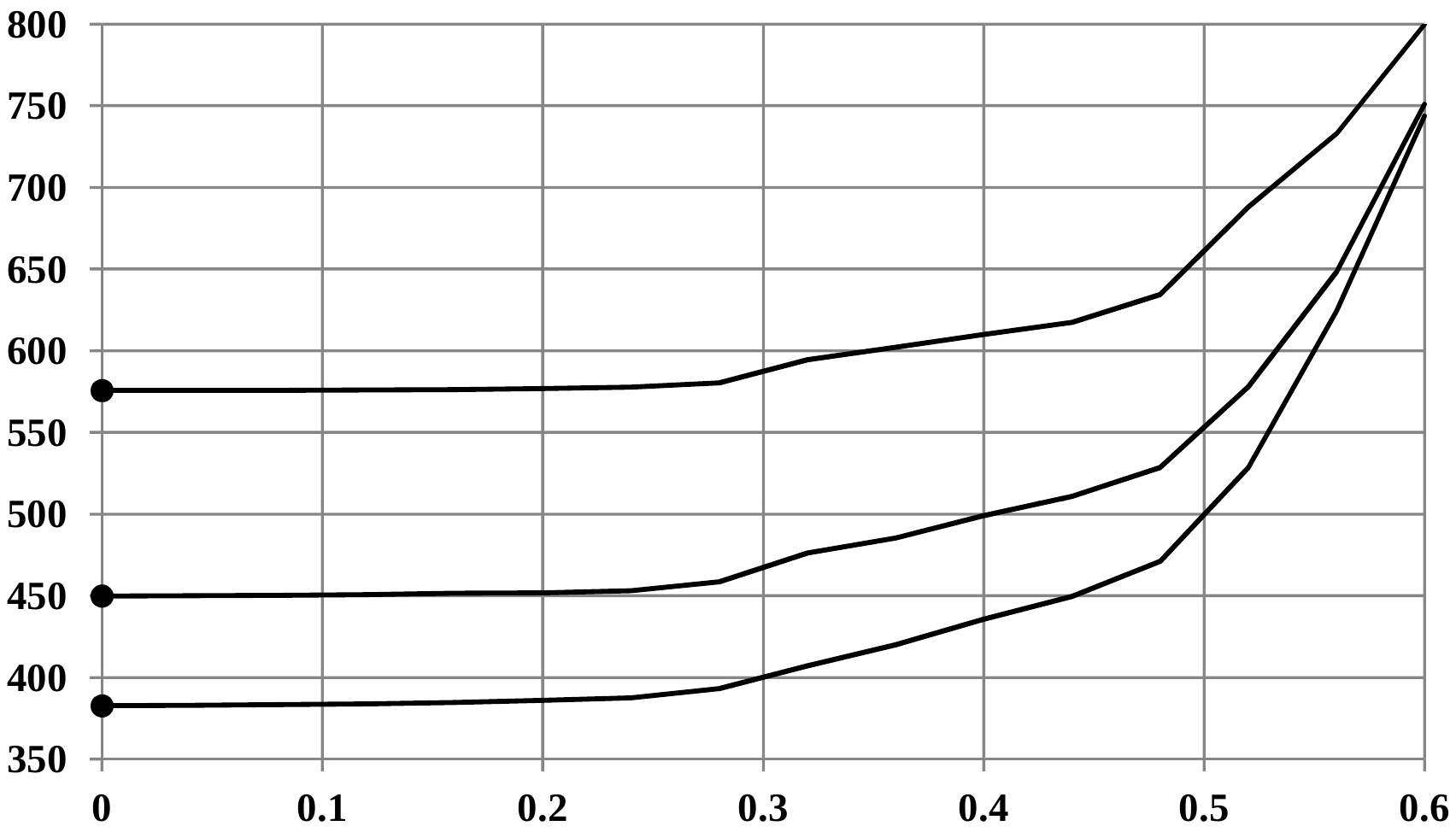}
		\scriptsize
		\tiny
		\put(-4.5,-0.1){$D$}
		\put(-1.5,-0.1){$D$}
		\put(-4.8,1.47){$p=2$}
		\put(-4.8,0.86){$p=3$}
		\put(-4.8,0.47){$p=4$}
		\put(-4.8,0.29){$p=5$}
		\put(-2.25,0.94){$p=10$}
		\put(-2.25,0.52){$p=15$}
		\put(-2.25,0.3){$p=20$}
		\end{picture}
		\caption{\label{ef500}The Efficient Frontiers for $n=500$}
	\end{center}
\end{figure}

\begin{figure}[ht!]
	\begin{center}
		\setlength{\unitlength}{1in}
		\begin{picture}(6,1.9)	
		\includegraphics[width=3in]{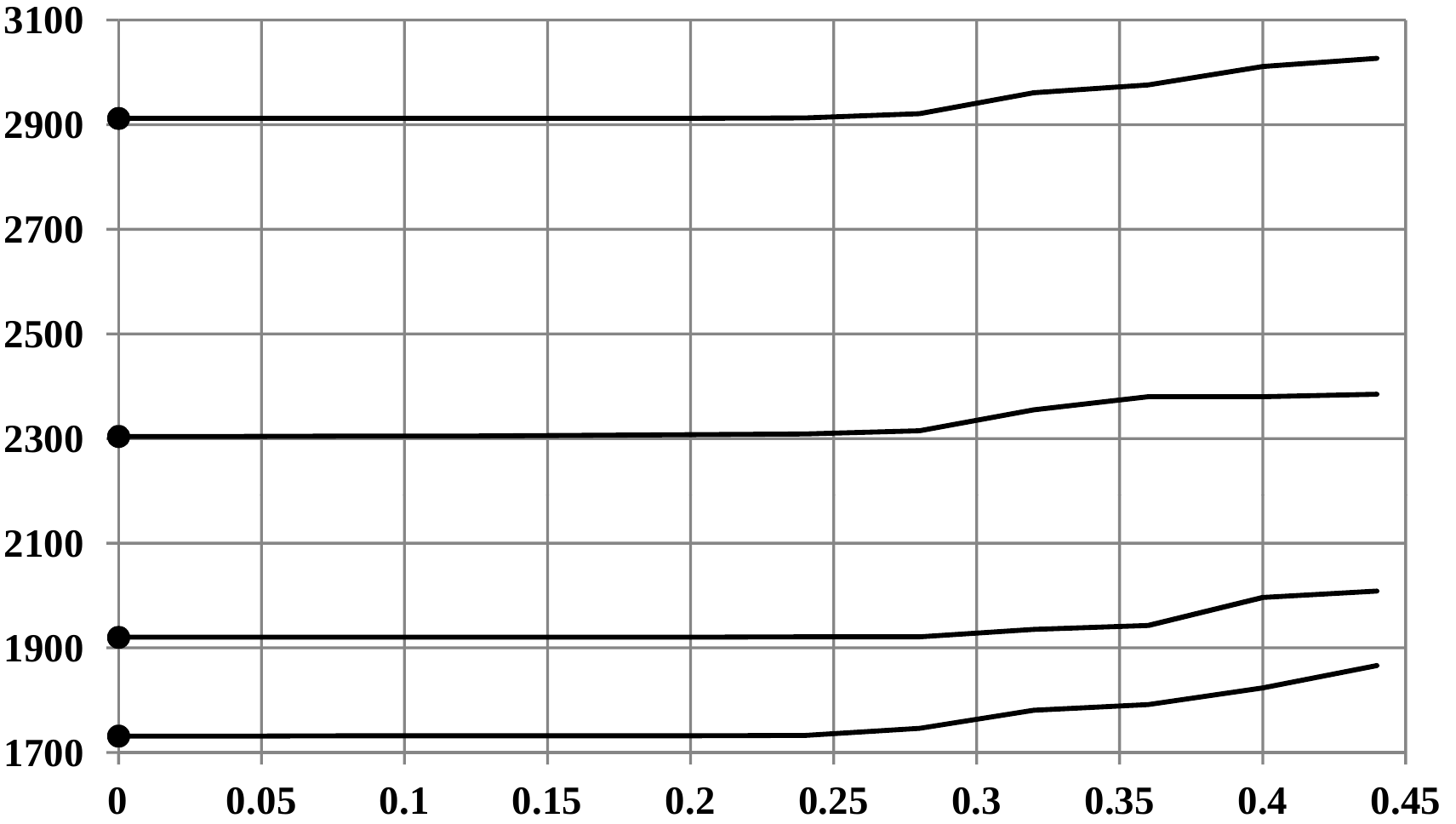}
		\includegraphics[width=3in]{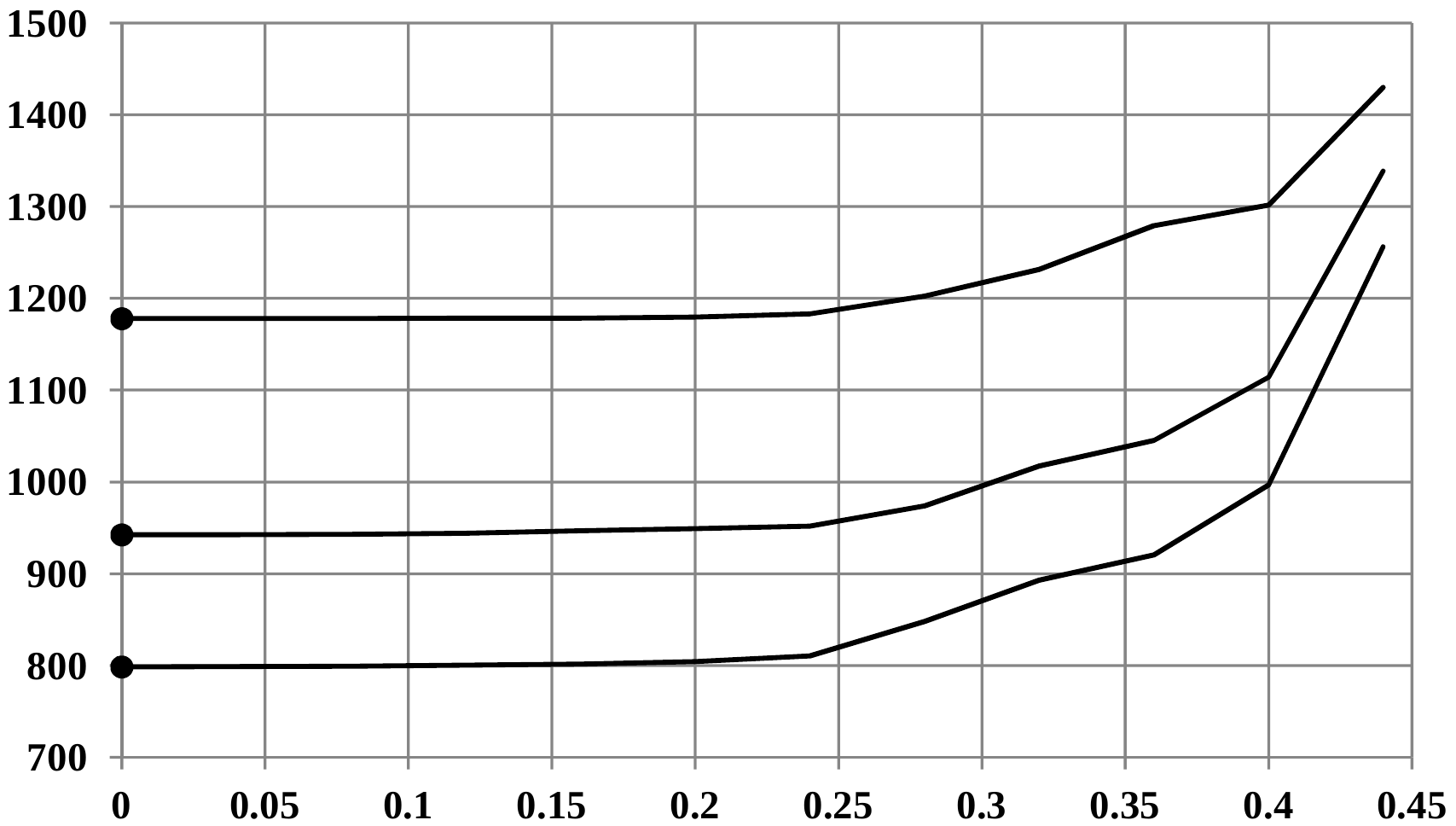}
		\scriptsize
		\tiny
		\put(-4.5,-0.1){$D$}
		\put(-1.5,-0.1){$D$}
		\put(-4.6,1.49){$p=2$}
		\put(-4.6,0.83){$p=3$}
		\put(-4.6,0.42){$p=4$}
		\put(-4.6,0.21){$p=5$}
		\put(-1.86,1.08){$p=10$}
		\put(-1.86,0.63){$p=15$}
		\put(-1.86,0.36){$p=20$}
		\end{picture}
		\caption{\label{ef1000}The Efficient Frontiers for $n=1000$}
	\end{center}
\end{figure}
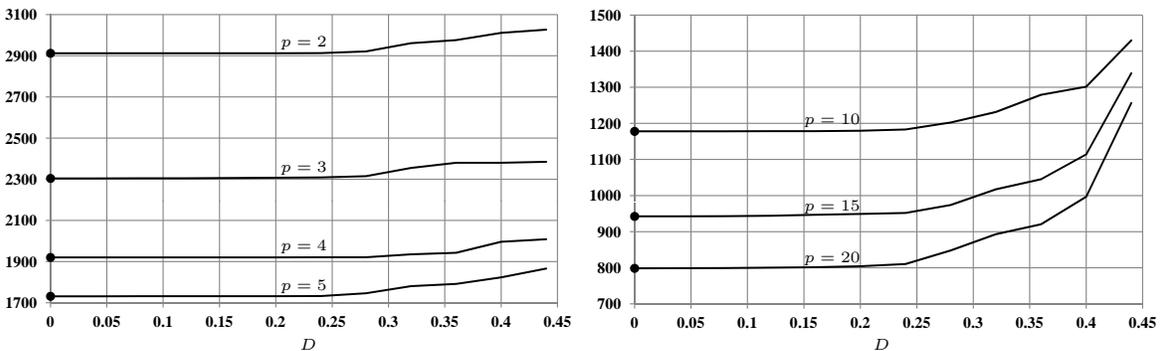

An interesting useful tool to support our model is the efficient frontier. For each $D$ we find the best solution and the efficient frontier will give the planner a trade-off between the minimum distance $D$ and the total cost. They may choose the best strategy by their own judgment. We constructed the efficient frontiers by applying SNOPT on the best Voronoi points. The efficient frontiers for $n=100, 500, 1000$ are depicted in Figures \ref{ef100}-\ref{ef1000}. As the number of facilities $p$ increases, the deterioration of the value of the objective function is faster when $D$ increases.

For example, in Figure \ref{ef100} the efficient frontiers for $p=2,3,4,5,10,15,20$, $n=100$ are depicted. Suppose that the plan is to build $p=15 $ facilities. When no minimum distance constraints are imposed ($D=0$, the black circle in the figure), the objective function is 74.47 (see Table \ref{Res}). When a minimum distance of $D=0.2$ is required, the objective function increases a bit to 75.04 which can be easily justified. When the planner wishes not to exceed a cost of 80 (see Figure \ref{ef100}), a minimum distance $D=0.43$ is possible yielding an objective of 79.81. $D=0.44$ has an objective of 80.03. If the planner is forced not to have facilities closer than $D=1$ mile from the demand points, the objective function is almost doubled to 136.74, and for $D=1.1$ it increases to 173.46. The planner has to apply his best judgment and select the best trade-off option. 

\section{\label{concl}Conclusions}

In this paper we proposed the obnoxious facilities $p$-median problem. A minimum given distance $D$ between facilities and a partial set of demand points is required. The objective is to find $p$ locations for facilities that minimize the weighted sum of distances between demand points and their closest facility subject to the minimum distance requirement. The resulting problem is extremely non-convex and traditional non-linear solvers applied in a multi-start approach do not perform well. An efficient solution method based on Voronoi diagrams is proposed and tested. It significantly outperformed standard non-linear solvers both in run times and quality of the solutions.

We also constructed the efficient frontiers of the test problems to assist planners in making location decisions. As expected, when the minimum required distance between facilities and demand points increases, the objective function deteriorates. The planner can use the efficient frontier in order to select the best trade-off solution.

\subsection{Suggestions for Future Research}

The obnoxious $p$-median problem can be formulated for other distance measures such as Manhattan ($\ell_1$) distances, general $\ell_p$ distances or location on the sphere. The methodology developed in this paper can be applied to such problems as long as the Voronoi points for such metrics can be obtained.

In a network environment we can replace the Voronoi points by the centers of links connecting demand points negatively affected by the facility that are at least $2D$ long. If such a link is less than $2D$ long, the whole link is infeasible. 

Locating obnoxious facilities minimizing other objectives can also apply the Voronoi points approach. For example, the obnoxious $p$-center problem is defined in a similar fashion. In the standard $p$-center problem \citep{KH79cen,CC09,CLY15}, rather than minimizing the sum of the weighted distances to the closest facility, the objective is minimizing the maximum weighted distance to the closest facility. The BLP formulation is changed to minimizing $L$ and adding  constraints $w_{i}d_{ij}y_{ij}\le L$. The MBLP formulation cannot be used because we must have $y_{ij}\in\{0,1\}$.

 The conditional version of the problem \citep{Mi80,BS90,OZ02} can be useful in many circumstances. This means that several facilities already exist in the area and $p$ additional facilities need to be located. Each demand point receives its services from the closest facility which may already exist or is a new one. The objective is  to minimize total cost subject to distance constraints. The BLP formulation is modified so that the existing facilities are assigned $x_j=1$ and are no longer variables, and the new facilities are restricted to the Voronoi points. An interesting variant is that some of the existing facilities violate the distance constraints and should be considered for removal.
 
		\renewcommand{\baselinestretch}{1}
		\renewcommand{\arraystretch}{1}
		\large
		\normalsize
		
		\bibliographystyle{apalike}
		

\end{document}